\documentclass[12pt]{article}
\usepackage{setspace}
\usepackage{textcomp}
\usepackage[dvips]{color} 
\usepackage{amsmath, amsthm, amssymb}
\usepackage[left=1.5in,top=.6in,right=1in,nohead]{geometry}
\usepackage{hyperref}
\usepackage{latexsym}
\usepackage{mathrsfs}
\usepackage{rawfonts}
\usepackage[dvips]{graphicx}
\usepackage[all]{xy}
\usepackage[sc,osf]{mathpazo}  

\def\bct{\begin{center}}
\def\ect{\end{center}}
\def\beg{\begin}

\def\<{\langle}
\def\>{\rangle}

\def\mbb{\mathbb}
\def\mbbr{\mathbb R}
\def\mbbz{\mathbb Z}

\def\ni{\noindent}

\def\tcb{\textcolor{blue}}
\def\tn{\textnormal}
\def\wt{\widetilde}

\newtheorem{thm}{Theorem}[section]

\def\A{{\mathbb A}}

\def\R{{\mathbb R}}
\def\L{{\mathbb L}}
\def\P{{\mathbb P}}

\def\Z{{\mathbb Z}}

\def\O{{\mathbb O}}
\def\H{{\mathbb H}}
\def\V{{\mathbb V}}
\def\Q{{\mathbb Q}}

\def\CP{\mathcal P}

\title{Algebraic topology of $G_2$ manifolds}
 \author{Selman Akbulut \and Mustafa Kalafat}

\begin{document}

\maketitle
\begin{abstract} In this paper we give a survey of various results about the topology of oriented Grassmannian bundles related to the exceptional Lie group $G_2$. 
Some of these results are new. We give self-contained proofs here. 
One often encounters these spaces when studying submanifolds of manifolds with calibrated geometries. 
For the sake of completeness we decided to collect them here in a self-contained
way to be easily accessible for future usage in calibrated geometry.
As an application we deduce existence of certain special 3 and 4 dimensional submanifolds of $G_2$ manifolds with special properties, which appear in the first named author's work with S. Salur about $G_2$ dualities.
\end{abstract}


\section{Introduction}

 Recall  that $G_2\subset SO(7)$ is the $14$-dimensional exceptional Lie group defined as the automorphisms of the imaginary octonions $im(\O)=\R^7$ preserving the cross product operation $\R^7 \times \R^7\to \R^7$ (e.g. \cite{hl}, \cite{bryant}, \cite{as1}, \cite{as2}). 
Octonions are the elements of the $8$ dimensional division algebra $\O=\H\oplus l \H= \R^8$  where $\H$ are the quaternions, $\O$ is generated by  $\langle 1, i, j, k, l, li ,lj, lk\rangle $. The cross product operation  $\times $ on $im (\O)$ is induced from the octonion multiplication on $\O$ by $u\times v=im(\bar{v} .u)$. 
We say an oriented $7$-manifold $M^7$ has a {\it $G_{2}$ structure} if its $SO(7)$- tangent frame bundle lifts to a $G_{2}$-bundle by the canonical fibration:  
$$G_{2}\to SO(7)\to\R\P^{7}\to BG_{2}\to BSO(7).$$
Alternatively $G_2$ can be defined by the special $3$-frames in $\R^7$ as follows:
\begin{equation*}
G_{2} =\{ (u_{1},u_{2},u_{3})\in (\R^{7})^3 \;|\; \langle u_{i},u_{j} \rangle=\delta_{ij}, \; \langle
u_{1} \times u_{2},u_{3} \rangle =0 \;  \},
\end{equation*}
or as  linear automorphisms of $\R^7$  preserving a certain  $3$-form $\varphi_{0}\in \Omega^{7}(\R^{7})$ 
\begin{equation*} G_{2}=\{ A \in GL(7,\R) \; | \; A^{*} \varphi_{0} =\varphi_{0}\; \} 
\end{equation*}
\vspace{.01in}
\noindent where 
$\varphi _{0}
=e^{123}+e^{145}+e^{167}+e^{246}-e^{257}-e^{347}-e^{356}$ with
$e^{ijk}=dx^{i}\wedge dx^{j} \wedge dx^{k}$.

\vspace{.05in}

 By using this last definition, a $G_{2}$ structure on $M^7$ can be defined as a 3-form $\varphi \in \Omega^{3}(M^{7})$ such that  at each $p\in  M$
the pair $ (T_{p}(M), \varphi (p) )$ is (pointwise) isomorphic to $(T_{0}(\R^{7}), \varphi_{0})$. This condition is equivalent to reducing the tangent frame bundle of a (not necessarily oriented) $7$-manifold $M$ from $GL(7, \R)$ to $G_2$.

\vspace{.1in}

The form $\varphi $  induces an orientation $\mu \in \Omega^{7}(M)$ on $M$, a metric $g=\langle ,\rangle_{\varphi }$ by $\langle u,v \rangle=[ i_{u}(\varphi ) \wedge i_{v}(\varphi )\wedge \varphi  ]/6\mu$, and $\varphi$ also defines a cross product operation $TM\times TM \mapsto TM$: $(u,v)\mapsto u\times v = u\times_{\varphi } v $ by $\varphi (u,v,w)=\<u \times v,w\>$. 

\vspace{.1in}

A manifold with $G_{2}$ structure $(M^{7},\varphi)$  is called a {\it $G_{2}$  manifold} (or an integrable $G_2$ structure) if at each point $p\in M$ there is an open chart  $(U,p) \to (\R^{7},0)$ on which $\varphi $ equals to $\varphi_{0}$ up to second order term, i.e. on the image of the open set $U$ we can write $\varphi (x)=\varphi_{0} + O(|x|^2)$. The condition that $(M,\varphi )$ be a $G_2$ manifold is equivalent to $\varphi$ being parallel under the induced metric connection $\nabla^{\varphi}(\varphi)=0$, which turns out to be equivalent to the condition $d \varphi=d^{*} \varphi =0$. 

\vspace{.1in }

Let $G_{k}^{+}\R^{n}$ denote the Grassmannian manifold or oriented $k$-planes in $\R^n$. We call $L\in  G_3^+\R^7$
an {\it associative $3$-plane} if $\varphi|_L\equiv vol(L)$. A 3-dimensional submanifold $Y\subset (M, \varphi)$ is called
{\em associative} if $\varphi|_Y\equiv vol(Y)$.
An equivalent condition of a submanifold $Y^3$ to be associative is that $\chi|_Y\equiv 0$,  where $\chi =\chi_{\varphi} \in \Omega^{3}(M,TM)$ is the  tangent bundle valued 3-form defined by  $ \langle \chi (u,v,w) , z \rangle=*\varphi  (u,v,w,z) $. This last identity implies a very useful property: $\chi$ assigns to every $3$-plane $L\subset TM$ an orthogonal vector $\chi|_{L}\in L^{\perp}\subset TM$. We also have:
\begin{equation*}
\varphi  (u,v,w)^2 + |\chi (u,v,w)|^2= |u\wedge v\wedge w|^2.
\end{equation*}
\begin{equation*}
\chi(u,v,w)= -u\times (v\times w)-\langle u,v\rangle w +\langle u,w\rangle v
\end{equation*}

We call $L\in G^{+}_{3}(\R^7) $ an {\it Harvey-Lawson $3$-plane} (HL plane in short) if $\varphi|_L\equiv 0$. We call  $S\in G^{+}_{4}(\R^7) $ a {\it coassociative $4$-plane} if $\varphi|_S\equiv 0$.
A $4$-dimensional submanifold $X^{4}\subset (M,\varphi)$  is called  {\em coassociative } if $\varphi|_X=0$. 
A manifold pair $(X^{4},Y^{3})$ such that  $Y^{3}\subset X^{4}\subset (M^{7}, \varphi)$ is called a {\it Harvey-Lawson pair} if the $\varphi \equiv 0$ on the restriction of the fibers of  the normal bundle $\nu(X)|_{Y}$ of  $X\subset (M,\varphi )$. The Grassmannians 
$G^{+}_{3}(\R^7)$ and $G^{+}_{4}(\R^7)$ have the following natural submanifolds 
\begin{equation*}ASS_{0}=\{L\in G^{+}_{3}(\R^7)\;|\; \varphi|_{L}=0\}\end{equation*} 
\begin{equation*}ASS_{+}=\{L\in  G^{+}_{3}(\R^7)\;|\; \varphi|_{L}=vol(L)\}\end{equation*} 
\begin{equation*}ASS_{-}=\{L\in G^{+}_{3}(\R^7)\;|\; \varphi|_{L}=-vol(L)\}\end{equation*}
\begin{equation*}COASS =\{S\in G^{+}_{4}(\R^7)\;|\; \varphi|_{S}=0\}\end{equation*}

\ni When there is no danger of confusion, we will abbreviate $ASS_{+}$ by 
$ASS$. Note that there is a natural identification $ASS\approx COASS$ 
given by $L \mapsto L^{\perp}$, and also

\begin{equation*}ASS_{\pm} \approx G_{2}/SO(4) \end{equation*}
\begin{equation*}ASS_{0} \approx G_{2}/SO(3) \end{equation*}

\ni From these descriptions it follows that $ASS_{0}$ is a sphere bundle 
over $ASS_{\pm}$ 
$$S^{3}\to ASS_{0}\to ASS_{\pm}$$
These special Grassmann manifolds sit in  $G^{+}_{3}(\R^7)$ as level 
sets of the function $$ \Phi : G^{+}_{3}(\R^7) \to \R $$ given by $L=u\wedge v\wedge w 
\mapsto \varphi_{0}(u,v,w)$, where $\{u,v,w\}$ is an orthonormal basis of $L$. 

We have $\Phi^{-1}(0)=ASS_{0}$ and $\Phi^{-1}(\pm1)=ASS_{\pm }$, since 
$\varphi_{0}$ is calibrating $3$-form 
$|\varphi_{0}(L)| \leq 1$. In this way $G^{+}_{3}(\R^7)$ appears as the 
double of the $D^{4}$-disk bundle over $ASS$. 
$\Phi$ is a Bott-Morse function. So that $ASS_{\pm}$ becomes two critical submanifolds 
with indices $0$ and $4$ respectively (\cite{jianweizhou}). \\

\begin{figure}[h] \bct \label{karpuz}
\includegraphics[width=.4\textwidth]{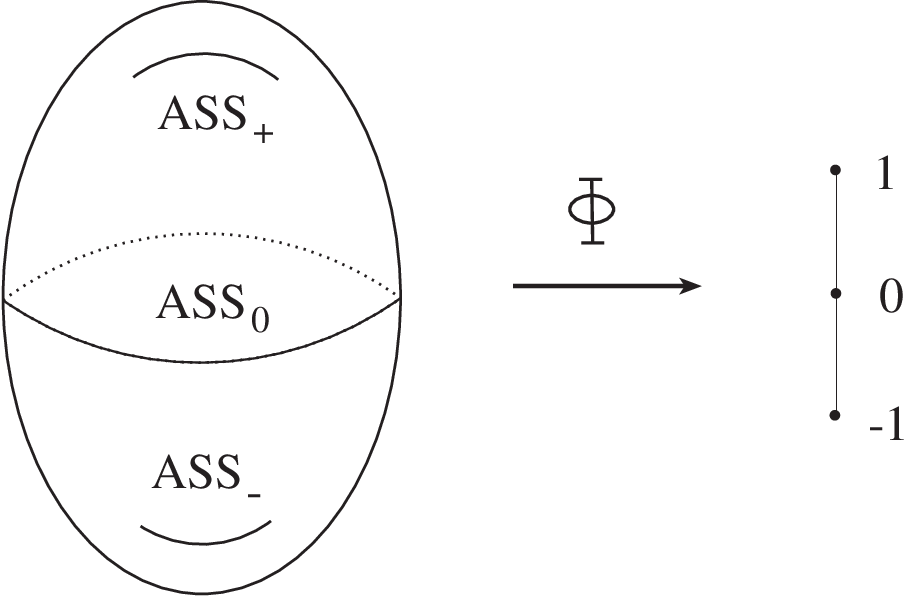}\\ \ect
  \caption{\, The map\, $\Phi: G^{+}_{3}(\R^7) \to \R$ } 
\end{figure}

\vspace{.1in}

These Grassmannians occur as the fibers of some bundles over $7$-manifolds with $G_{2}$ structure $(M^{7},\varphi)$, providing a useful tool studying deformations of associative submanifolds \cite{as2}.  Next we summarize some of the constructions from \cite{as1}. 
Let  $\P_{SO(7)}\to M$ be frame bundle of the tangent bundle $T(M)\to M$ of any closed smooth oriented $7$-manifold $M$, and let $\widetilde{M}\to M$ be the bundle oriented $3$-planes in $TM$, which
is defined by the identification  $[p,L]=[pg,g^{-1}L] \in
\widetilde{M}$:

\begin{equation*}
 G^{+}_{3}(\R^7) \to \widetilde{M}= \P_{SO(7)}(M)\times_{SO(7)} G^{+}_{3}(\R^7) \overset{\pi}{\longrightarrow} M.
\end{equation*}

\ni Let $\xi \to G^{+}_{3}(\R^7)$, and $\nu =\xi^{\perp}\to G^{+}_{4}(\R^7)$ be the universal $\R^3$ bundle, and its dual $\R^{4}$ bundle, respectively. Therefore, $Hom(\xi,
\nu)=\xi^{*}\otimes \nu \longrightarrow G^{+}_{3}(\R^7)$ is the tangent
bundle  $TG^{+}_{3}(\R^7)$. $\xi$, $\nu$ extend fiberwise to give bundles
$\Xi \to  \widetilde{M}$, $\V \to \widetilde{M}$ respectively. If
$\Xi^*$ be the dual of $\Xi$,  then $ Hom (\Xi, \V)=
\Xi^{*}\otimes \V \to \widetilde{M}\; $ is the bundle  of vertical
vectors  $T^{v}(\widetilde{M}) $ of $T(\widetilde{M}) \to M$, i.e. the
tangents to the fibers of  $\pi $.

\vspace{.1in} 

When $(M^7, \varphi)$  is a manifold with $G_2$ structure, similarly to the construction above,   we can form  the following subbundles of $\widetilde{M}\to M$ 
\begin{equation*}
ASS_{\pm}\to \A SS_{\pm} = \P_{G_{2}}(M) \times _{G_{2}}ASS_{\pm} \longrightarrow M
\end{equation*}
\begin{equation*}
ASS_{0}\to \A SS_{0} = \P_{G_{2}}(M) \times _{G_{2}}ASS_{0} \longrightarrow  M
\end{equation*}
where $ \P_{G_{2}}(M)$ is the $G_2$ frame bundle of the tangent bundle of $(M,\varphi)$. In particular  $\A SS=\P_{G_{2}}(M)/SO(4)\longrightarrow \P_{G_{2}}(M)/G_{2} =M$. As in the previous section,  the restriction  of the universal bundles $\xi, \;\nu=\xi^{\perp} \to G^{+}_{3}(\R^7)$ induce $3$  and $4$ plane bundles $\Xi\to\A SS$ and $\V\to \A SS$.   Also we have the similar map $L\to \varphi(L)$
\begin{equation*}
\Phi : \widetilde{M}\longrightarrow \R
\end{equation*}
with $\Phi^{-1}(0)=\A SS_{0}$ and $\Phi^{-1}(\pm1)=\A
SS_{\pm }$. Fiberwise this is just the map previously described on $G^{+}_{3}(\R^7)$, it is the bundle version of the map described in Figure~\ref{karpuz}. So we have disjointly embedded pair of codimension $4$-submanifolds $\A SS_{\pm} \subset \widetilde{M}$, which are separated by a codimension zero submanifold $\A SS_{0}\subset \widetilde{M}$.

\vspace{.1in}

Any embedding of a $3$-manifold  $f: Y^{3}\hookrightarrow M^7 $, by its tangential Gauss map,  
lifts to an embedding $f_{T}: Y \hookrightarrow \widetilde{M}$ such that the pull-backs $f_{T}^{*} \Xi=T (Y)$ and  $f_{T}^{*} \V=\nu(Y)$ are the tangent and normal bundles of $Y$.  In particular,   if $f$ is an embedding of an  associative submanifold ,  then the image of $f_{T}$ lands in $\A SS $ 
\begin{equation*}
\begin{array}{lcl}
  & & \widetilde{M}\,\supset\, \A SS_{0},\; \A SS_{\pm} \\
  \hspace{.25in} f_{T}   \hspace{-.1in}&  \nearrow \;  &   \downarrow\\
\;Y\;\;& \stackrel{f}{\longrightarrow}   & M
\end{array}
\end{equation*}
Similarly any embedding of a $4$-manifold  $f: X^{4}\hookrightarrow M^7 $, by its normal Gauss map, induces an embedding $f_{N}: X\hookrightarrow \widetilde{M}$, such that
$f_{N}^{*} \Xi=\nu(X)$ and  $f_{N}^{*} \V=T(X)$ are the normal and tangent bundles of $X^4$.

\vspace{.1in}

 If $\L=\Lambda^{3}(\Xi)\to \widetilde{M}$  is the determinant (real) line  bundle. By the discussion above $\chi $ maps every oriented $3$-plane in $T_{x}(M)$ to its $4$-dimensional complementary subspace,  so $\chi$ gives a bundle map $\L \to \V$ over $\widetilde{M}$, which is a section of $\L^{*}\otimes \V \to \widetilde{M}.$  Since $ \Xi $ is oriented  $\L$ is trivial,   so $ \chi $ actually gives  a section

\begin{equation*}
 \chi =\chi_{\varphi} \in \Omega^{0}(\widetilde{M},  \V)
 \end{equation*}

\noindent $\A SS$  is the zeros of this section.  Associative submanifolds $Y\subset M$ are characterized by the condition $\chi|_{\tilde{Y}}=0$,  where $\tilde{Y}\subset \widetilde{M}$ is the canonical lifting of $Y$.

\vspace{.1in}

$\A SS $ is the universal space
parameterizing associative submanifolds of $M$. 
In particular,  if $\tilde{f}: Y\hookrightarrow  \widetilde{M}_{\varphi}$ is the lifting of an associative submanifold, by pulling back we see that  the principal  $SO(4)$   bundle $\CP(\V)\to \widetilde{M}_{\varphi} $ induces an $SO(4)$-bundle $\CP (Y)\to Y$, and gives the following vector bundles  via the  representations:

\begin{equation}
\begin{array}{lcc}
\; \nu(Y) \hspace{.2in}:   &y \mapsto  q y\lambda^{-1}   & \hspace{1in}  \\
 \; T (Y)  \hspace{.1in}:  & x \mapsto  qx q^{-1} &   \\
 \end{array}
 \end{equation}
where $[q,\lambda]\in SO(4)=SU(2)\times SU(2)/\Z_{2}$, $\nu=\nu(Y)$ and 
$T(Y)= \lambda_{+}(\nu ) $. Here we use quaternionic notation $\R^4 =\Q$ 
and $\R^3= im(\Q)$.  Also we can identify $T^{*}Y$ with $TY$ by the 
induced canonical metric coming from $G_2$ structure. From above we  have 
the action $T^{*}Y \otimes \nu \to \nu  $   inducing  actions 
$\Lambda^{*}(T^{*}Y) \otimes \nu \to \nu  $.

\vspace{.1in}

In the next section we survey relevant results from algebraic topology of 
these spaces. Most of these results are elementary, folklore or already known (e.g. \cite{borel}, \cite{mimura}, \cite{zhou}), but some results and the application in Section \S\ref{hlpair} 
are new. There are also other applications e.g. \cite{ku}.   
For the sake of completeness we decided to collect them here in a self contained way to be easily accessible for future usage in calibrated geometry. 
 
 \vspace{.1in}

{\bf Acknowledgements.} This work was partially supported by NSF (National Science Foundation) grant DMS-1065879, FRG-1065827 and T\"ubitak (Turkish science and research council) grant {$\sharp$}114F320. Thanks to \.I. \"Unal for some useful discussions. 


\section{Algebraic topology of Grassmann bundles}

Here compute cohomology groups of various Grassmann bundles. To this end we first start with the following calculations.  

\vspace{.1in}

\noindent {\bf Lemma \ref{g27homology}.}
{\em The homology of the Grassmann manifold $G_2^+\mbb R^7$ of oriented $2$-planes is given by the following table: 
$$H_*(G_2^+\mbb R^7;\mbb Z)=(\mbbz,0,\mbbz,0,\mbbz,0,\mbbz,0,\mbbz,0,\mbbz).$$}

 This result is obtained by some elementary computations on various spectral sequences. After this, by computing the torsion and using the Gysin sequence of some special fibrations we get the following.
 
\vspace{.1in}

\noindent {\bf Theorem \ref{g37homology}.} 
{\em The homology of the oriented Grassmann manifold $G_3^+\mbb R^7$ is given by: 
$$H_*(G_3^+\mbb R^7;\mbb Z)=
(\mbbz,0,\mbbz_2,0,\mbbz\oplus\mbbz,\mbbz_2,\mbbz_2,
0,\mbbz\oplus\mbbz,\mbbz_2,0,0,\mbbz).$$ }

\noindent Also the cup product structure of the cohomology can be computed as follows. 

\vspace{.1in}

\noindent {\bf Theorem \ref{cupproductstructuregrassmannian}.} 
{\em The cohomology ring of  $G_3^+\mbbr^7$ is given by: 
$$H^*(G_3^+\mbbr^7;\mbbz)=\mbbz[x,y]/(x^4,x^3-y^3,x^2-y^2,xy-yx) 
\oplus\mbbz_2[z,t]/(z^3,t^2,zt+tz)$$ 
where the degrees of the generators are~ $|x|=|y|=4, |z|=3$ ~and~ $|t|=7$.}

\vspace{.1in}

\noindent  As a by product, along the way we compute the homology of a Stiefel manifold.

\vspace{.1in}

\noindent {\bf Theorem \ref{v37homology}.} 
{\em The homology of the Stiefel manifold $V_3\mbb R^7$ is given by: 
$$H_*(V_3\mbb R^7;\mbb Z)=
(\mbbz,0,0,0,\mbbz,\mbbz_2,0,0,0,\mbbz_2,0,\mbbz,
0,0,0,\mbbz).$$ }


\noindent Then by combining two different fibrations we can compute the homology of the Lie group $G_2$, and compute its cohomology ring as well. 

\vspace{.05in}

\ni {\bf Theorem \ref{g2homology}.} 
{\em The homology groups of the Lie group $G_2$ are given as follows:
$$H_*(G_2;\mbb Z)=
(\mbbz,0,0,\mbbz,0,\mbbz_2,0,0,\mbbz_2,0,0,\mbbz,0,0,\mbbz).$$}


\vspace{.05in}

\ni {\bf Theorem \ref{g2cup}.} 
{\em The cohomology ring of the Lie group $G_2$ can be described as follows.
$$H^*(G_2;\mbb Z)=\Lambda_\mbbz[x_3,x_{11}]\oplus
\Lambda_{\mbbz_2}[x_6,x_9]/(x_6 x_9)$$
where the degrees of the generators are~ $|x_k|=k$.
}

\vspace{.1in}

In part \ref{sectiongrassmanninvariants} we deal with topological 
computations on the Grassmannian, in part \ref{sectioncupproduct} we 
compute the cup product structure. 
Finally, in part \ref{sectiong2homology} we analyze  the space $G_2$  
and compute its cohomology ring, and then in \ref{sectionclassifying} compute cohomology rings of of certain bundles associated to $G_2$ manifolds.
\vspace{.05in}

\subsection{Homology of Grassmannians}\label{sectiongrassmanninvariants}
The aim of this section is to compute some integral homology groups 
of the oriented real Grassmannian 
$G_3^+\mbb R^7$. 
We will use various forms of the Serre spectral sequence of a fiber bundle with various coefficients. 
As a warm up, let us recall the homology and cohomology of the basic spaces. Starting from $SO_3$, using its identification with $\mbb{RP}^3$ and a combination of Poincar\'e duality with the universal coefficients theorem (UCT) one easily computes its homology groups as
$$H_*(SO_3;\mbbz )  = ( \mbbz   , \mbbz_2 ,  0     , \mbbz   ).$$
$$H_*(SO_3;\mbbz_2) = ( \mbbz_2 , \mbbz_2 ,\mbbz_2 , \mbbz_2 ).$$
Now, consider the Stiefel manifold $V_3\mbb R^7$ which is defined to be the 
space of orthonormal 3-frames in $\mbb R^7.$ By the Stiefel fibrations (e.g. \cite{hatcher}) this is a $7\negthickspace -\negthickspace3\negthickspace-\negthickspace1\negthickspace=\negthickspace3$ connected space. Since the dimension is even the fourth homotopy group is $\mbbz$ by Stiefel \cite{stiefel} so $\pi_{01234}(V_3\mbb R^7)=(0,0,0,0,\mbbz )$. 
This notation expresses the homotopy groups of the space up to level 4. 
See \cite{whitehead,paechter} for higher homotopy. Sending a 3-frame to the oriented 3-plane, which
it spans, gives us the fibration
\beg{equation} \label{fibervg}
SO_3 \longrightarrow V_3\mbb R^7 \longrightarrow G_3^+\mbb R^7.\end{equation} 
Using the related homotopy exact sequence (HES) and homotopy groups
$$\pi_{01234}(SO_3)=( 0,\mbbz_2,0,\mbbz,\mbbz_2 )$$
$$\pi_{01234}(G_3^+\mbb R^7)=(0,0,\mbbz_2,0,\mbbz\oplus\mbbz).$$
Our next aim is to compute some homology groups for the Grassmannian. 
From above by the Hurewicz isomorphisms
$$H_{012}(G_3^+\mbb R^7)=(\mbbz,0,\mbbz_2).$$
The Poincar\'e polynomial of $G_3^+\mbb R^7$ is known to be 
\beg{equation} \label{poincarepolynomialg3r7}
p_{G_3^+\mbb R^7}(t) = 1 + 2t^4 + 2t^8 + t^{12}.
\end{equation}
See \cite{ghv} vol.III, pp.494-496 for computations also \cite{gluckmackenziemorgan}. For further
homology computations, we will consult to the spectral sequences 
and the Gysin sequence. 
We will abbreviate $G=G_3^+\mbb R^7$ and $V=V_3\mbb R^7$ frequently in what follows. 

\noindent Our first assertion is the following Lemma.
\beg{lem}\label{H_3G}For the oriented Grassmann manifold we have $H_3(G_3^+\mbb R^7;\mbb Z)=0$.
\end{lem}
\beg{proof} 
We consider the homological Serre spectral sequence with $\mbb Z$-coefficients
associated to the fiber bundle \eqref{fibervg}, properties of which is given as follows \cite{hajimesato}.
$$ E_{p,q}^2 := H_p( G ; H_q(SO_3;\mbb Z) ) $$
$$ E_{p,q}^\infty = { F_{p,q} / F_{p-1,q+1} }$$
where $F_{p,q}$ are abelian groups forming a filtration satisfying
$$0 = F_{-1,n+1} \subset\cdots\subset F_{n-1,1} \subset F_{n,0} = H_n(V;\mbb Z).$$
The differentials are bidegree $(-n,n-1)$ maps 
$$d^n : E_{p,q}^n \longrightarrow E_{p-n,q+(n-1)}^n.$$
\beg{table}[ht] 
\label{table:E2Hom}
\caption{ {\em Homological Serre spectral sequence for $G_3^+\mbb R^7$, 
second page.}}
\bct{
$\begin{array}{cc|c|c|c|@{} c @{}|c|c|} \cline{3-8}
    &3&       &       &   &&&\\ \cline{3-8}
    &2& \,0\, &       &   &&&\\ \cline{3-8}
    &1&       & \,0\, & ~~ &~~&~~&\\ \cline{3-8}
E^2~~~&0&     &       &   & H_3(G;\mbbz) &&  \\ \cline{3-8}
& \multicolumn{1}{c}{ } & \multicolumn{1}{c}{0} 
& \multicolumn{1}{c}{1} & \multicolumn{1}{c}{2} 
& \multicolumn{1}{c}{3} & \multicolumn{1}{c}{4}
& \multicolumn{1}{c}{5}   \\  
\end{array}$   
\setlength{\unitlength}{1mm}
\begin{picture}(0,0)
\put(-31.5,-2){\vector(-3,1){8}}   
\end{picture} }
\ect
\end{table}
\noindent Some of the terms appear 
in the following table.
Note that 
$$E^2_{1,1}=H_1(G;H_1(SO_3;\mbbz))=H_1(G;\mbbz_2)=H_1\otimes\mbbz_2\oplus\tn{Tor}(H_0;\mbbz_2)=0.$$
We have immediate convergence for the $E^2_{3,0}$ term so that
$$H_3(G;\mbbz)=E^2_{3,0}=\cdots=E^\infty_{3,0}=F_{3,0}/F_{2,1}=H_3(V;\mbbz)/F_{2,1}=0.$$ \end{proof}

\noindent Another result on the Grassmannian is the following.

\beg{lem}\label{g27homology}The homology of the oriented Grassmann manifold $G_2^+\mbb R^7$ is given by: 
$$H_*(G_2^+\mbb R^7;\mbb Z)=(\mbbz,0,\mbbz,0,\mbbz,0,\mbbz,0,\mbbz,0,\mbbz).$$
\end{lem}

\vspace{.05in}

\beg{proof} Consider the fibration
\beg{equation} S^1=SO_2\longrightarrow V_2\mbb R^7 \longrightarrow G_2^+\mbb R^7.
\label{fibervg27} \end{equation}
The homology of the Stiefel manifold is well-known (e.g. \cite{hatcher}), it is given by: 
\beg{equation}\label{v27homology}
H_*(V_2\mbb R^7;\mbb Z)=(\mbbz,0,0,0,0,\mbbz_2,0,0,0,0,0,\mbbz)
\end{equation}
Since $V_2\mbb R^7$ is $7-2-1=4$ connected, homotopy exact sequence of the above fibration immediately gives 
$\pi_{012}(G_2^+\mbb R^7)=(0,0,\mbb Z)$.    
The homological Serre spectral sequence reads as follows.
$$ E_{p,q}^2 := H_p( G_2^+\mbb R^7 ; H_q(S^1;\mbb Z) ) $$
$$ E_{p,q}^\infty = { F_{p,q} / F_{p-1,q+1} }$$
where $F_{p,q}$ are abelian groups of the filtration 
$$0 = F_{-1,n+1} \subset\cdots\subset F_{n-1,1} \subset F_{n,0} = H_n(V_2\mbb R^7;\mbb Z).$$

\noindent We first fill out the limiting page of the sequence as in the Table \ref{table:EinftyHomG27}, except the terms in quotation marks, which we do not use, however the reader can compute them after secondary steps and written here for recording purposes only. 
Vanishing of the homology of $V_2\mbb R^7$ and the filtrations easily handle this far.
\beg{table}[ht]
\caption{ {\em The limiting page of Homological Serre spectral sequence for $G_2^+\mbb R^7$.}}
\bct{
$\begin{array}{cc|c|c|c|c|@{}c@{}|@{}c@{}|c|c|c|c|c|} \cline{3-13}
           &1&0&0&0&0&\tn `\mbbz_2\hspace{-.8mm}\tn{'}&0&0&0&0&0& \\ \cline{3-13}
E^\infty~~~&0& &0&0&0&   0   &`0\tn'&0&0&0&0&0  \\ \cline{3-13}
& \multicolumn{1}{c}{ } & \multicolumn{1}{c}{0} 
& \multicolumn{1}{c}{1} & \multicolumn{1}{c}{2} 
& \multicolumn{1}{c}{3} & \multicolumn{1}{c}{4}
& \multicolumn{1}{c}{5} & \multicolumn{1}{c}{6} 
& \multicolumn{1}{c}{7} & \multicolumn{1}{c}{8}
& \multicolumn{1}{c}{9} & \multicolumn{1}{@{}c@{}}{10} \\  
\end{array}$   
}\ect \label{table:EinftyHomG27}  \end{table}

\noindent Next we fill out the second page. Keep in mind that in the following the entries of a column are identical. Columns till the second one follows from the homotopy groups. The third column is zero since we have the immediate
convergence $E^2_{3,0}=E^\infty_{3,0}$ because of the differential.
The fourth column is the outcome of the isomorphism
$$ 
\mbbz~ \tilde \longleftarrow ~E_{4,0}^2  ~:~  d^2_{4,0}$$
Since the domain and image group have to converge to zero in the next page, this map is both injective and surjective. The sixth column and on are the consequences of the universal coefficients theorem and the Poincar\'e duality. 
Fifth column is the remaining one. Again, because of the immediate convergence
$E_{5,1}^2=E_{5,1}^\infty$ the $E_{5,1}^2$ term has to vanish, 
\beg{table}[ht]
\caption{ {\em 
The second page.}}
\bct{
$\begin{array}{cc|c|c|c|c|c|c|c|c|c|c|c|} \cline{3-13}
    &1&\mbbz&0&\mbbz&0&\mbbz&0&\mbbz&0&\mbbz&0&\mbbz \\ \cline{3-13}
E^2~~~ &0&  &0&\mbbz&0&\mbbz& &\mbbz&0&\mbbz&0&\mbbz \\ \cline{3-13}
& \multicolumn{1}{c}{ } & \multicolumn{1}{c}{0} 
& \multicolumn{1}{c}{1} & \multicolumn{1}{c}{2} 
& \multicolumn{1}{c}{3} & \multicolumn{1}{c}{4}
& \multicolumn{1}{c}{5} & \multicolumn{1}{c}{6} 
& \multicolumn{1}{c}{7} & \multicolumn{1}{c}{8}
& \multicolumn{1}{c}{9} & \multicolumn{1}{@{}c@{}}{10} \\  
\end{array}$   
\setlength{\unitlength}{1mm}
\begin{picture}(20,10)
\put(-59,2){\vector(-2,1){8}}
\end{picture} }
\ect\label{table:E2HomG27}\end{table}
so is this column. A row of Table \ref{table:E2HomG27} determines the homology of $G_2^+\mbb R^7$ by its definition.
\end{proof}

Now we can prove our main Lemma.
\beg{lem} The torsion subgroup of $H_4(G_3^+\mbb R^7;\mbb Z)$ is trivial.
\end{lem}
\beg{proof}
Since the free part $F_5$ of $H_5(G;\mbbz)$ is trivial \ref{poincarepolynomialg3r7}, we can compute  the
torsion part of $H_4(G;\mbbz)$ which is denoted by $T_4$ using cohomology as follows.
$$H^5(G
; \mbbz)=\tn{Hom}(H_5,\mbbz)\oplus \tn{Ext}(H_4,\mbbz)=0\oplus T_4=T_4.$$
To figure out this group we will work with two new fibrations.
$$\xymatrix{
           & S^4\ar[d]& \\
S^2 \ar[r] & S(E_3\mbb R^7)=E_0 \ar[d]^{\pi_v}\ar[r]^{~~~~~~\pi_h} &G_3^+\mbb R^7\\
  &G_2^+\mbb R^7 &  }$$ 
\noindent Here $E_3\mbb R^7$ denotes the tautological bundle over 
$G_3^+\mbb R^7$. The horizontal fibration is clear, which is obtained by 
removing the zero section. From there one can obtain the vertical 
fibration
with the following procedure. A point in $G_2^+\mbb R^7$ represents a 2-
plane which is contained in $(7-2)(3-2)=4$   parameter of 3-planes. Since 
we take the orientations into consideration we obtain spheres rather than 
projective spaces.  
One may think in terms of the oriented flag variety
$F_{2,3}^+(\mbb R^7)$ with its projection maps. See \cite{harris}. 
Now, consider the Gysin exact sequence \cite{milnor} of the vertical 
fibration.
$$\cdots\longrightarrow
H^0(G_2^+\mbb R^7;\mbbz)\stackrel{\cup e_5}\longrightarrow 
H^5(G_2^+\mbb R^7;\mbbz)  \stackrel{\pi_v^*} \longrightarrow 
H^5(E_0          ;\mbbz)   \longrightarrow 
H^1(G_2^+\mbb R^7;\mbbz) \longrightarrow \cdots 
$$
Since we proved that the odd homology of the grassmannian $G_2^+\mbb R^7$ 
is zero in Lemma \ref{g27homology}, this implies that the middle term 
$H^5(E_0;\mbbz)$ vanishes.

\noindent Next consider the Gysin sequence of the horizontal fibration.
$$\cdots\longrightarrow
H^2(G_3^+\mbb R^7;\mbbz)\stackrel{\cup e_3}\longrightarrow 
H^5(G_3^+\mbb R^7;\mbbz)  \stackrel{\pi_h^*} \longrightarrow 
H^5(E_0          ;\mbbz)   \longrightarrow 
\cdots 
$$
Since $H_2(G_3^+\mbb R^7;\mbbz)$ is torsion and $H_1(G_3^+\mbb R^7;\mbbz)$ 
is zero we have $H^2(G_3^+\mbb R^7;\mbbz)=0.$ Together with the vanishing 
of 
$H^5(E_0;\mbbz)$ we obtain our result
$$H^5(G_3^+\mbb R^7;\mbbz)=0.$$        \end{proof}

\noindent Since by \ref{poincarepolynomialg3r7} the fourth Betti number of 
$G_3^+\mbb R^7$ is $2$, this Lemma implies the following.
\beg{cor}\label{H_4G} For the Grassmann manifold we have $H_4(G_3^+\mbb 
R^7;\mbb Z)=\mbb Z \oplus \mbb Z$.
\end{cor}

The results so far helps us to consume most of the homology of our 
Grassmann manifold. The homology at the levels $7$ and above are easily 
deduced from the Ext universal coefficients theorem and Poincar\'e 
duality. Finally using (\ref{poincarepolynomialg3r7}) in addition yields 
that the homology in the levels $5$ and $6$ are solely torsion, isomorphic 
and denoted by $T_5$. The rest of this section is devoted to compute this 
group. 

In order to compute this torsion this time we need some results on some 
Stiefel manifolds.

\beg{thm}\label{v37homology}The homology of the Stiefel manifold $V_3\mbb 
R^7$ is computed as, 
$$H_*(V_3\mbb R^7;\mbb Z)=
(\mbbz,0,0,0,\mbbz,\mbbz_2,0,0,0,\mbbz_2,0,\mbbz,
0,0,0,\mbbz).$$
\end{thm}
\beg{proof} We will be using the homological Serre spectral sequence 
related to the following new fibration.
$$S^4\longrightarrow V_3\mbb R^7 \longrightarrow V_2\mbb R^7.$$
This is obtained by projecting onto the first two vectors of the frame 
and the third one has unit independency in $\mbb R^5$. Defining groups are 
as follows.  
$$ E_{p,q}^2 := H_p( V_2\mbb R^7 ; H_q S^4 ) $$
$$ E_{p,q}^\infty = { F_{p,q} / F_{p-1,q+1} }$$
where $F_{p,q}$ are abelian groups of the filtration 
$$0 = F_{-1,n+1} \subset\cdots\subset F_{n-1,1} \subset F_{n,0} = 
H_n(V_3\mbb R^7;\mbb Z).$$
Merely knowing the homology of $V_2\mbb R^7$ as in (\ref{v27homology})
one can construct the 
second page of the spectral sequence. 
\beg{table}[ht]
\caption{ {\em Homological Serre spectral sequence for $V_3\mbb R^7$ %
with $\mbb Z$-coefficients. 
}}
\bct{
$\begin{array}{cc|c|c|c|c|c| c |c|c|c|c|@{}c@{}| @{} c @{} |} \cline{3-14}
  &4& \mbbz &&&&& \mbbz_2 &&&&&& \mbbz\\ \cline{3-14}
  &3&       &&&&&         &&&&&&      \\ \cline{3-14}
  &2&       &&&&&         &&&&&&      \\ \cline{3-14}
  &1&       &&&&&         &&&&&&      \\ \cline{3-14}
E^2=E^\infty~~~
  &0& \mbbz &&&&& \mbbz_2 &&&&&& \mbbz\\ \cline{3-14}
& \multicolumn{1}{c}{ } & \multicolumn{1}{c}{0} 
& \multicolumn{1}{c}{1} & \multicolumn{1}{c}{2} 
& \multicolumn{1}{c}{3} & \multicolumn{1}{c}{4}
& \multicolumn{1}{c}{5} & \multicolumn{1}{c}{6} 
& \multicolumn{1}{c}{7} & \multicolumn{1}{c}{8} 
& \multicolumn{1}{c}{9} & \multicolumn{1}{c}{10}   
& \multicolumn{1}{c}{11}
\\  
\end{array}$   
\setlength{\unitlength}{1mm}
\begin{picture}(0,0)
\put(-47.5,-6){\vector(-2,1){9}}   
\end{picture}
}\ect
\end{table}
Because of the abundance of zeros and freeness of $\mbbz$ we have the 
immediate convergence. Isomorphism of the groups in the filtration gives 
the triviality of the homology of $V_3\mbb R^7$ at the levels 1 to 3,  
6 to 8, 10 and 12 to 14 since their diagonal consist entirely of zeros.    
Now, to determine the homology at the 4th level start at $F_{-1,5}=0$. 
Since $E^\infty_{0,4}=F_{0,4}/F_{-1,5}=\mbbz$ implying that 
$F_{0,4}\approx\mbbz$. 
Next 
$$F_{1,3}/F_{0,4}=E^\infty_{1,3}=\cdots =E^\infty_{4,0}=
F_{4,0}/F_{3,1}=0$$ implying the isomorphisms 
$$\mbbz \approx F_{0,4}\approx F_{1,3}\approx F_{2,2}\approx 
F_{3,1}\approx F_{4,0}=H_4(V_3\mbb R^7;\mbb Z).$$
In an exactly similar way the other two nontrivial groups on the 4th row 
projects onto the homology at levels 9 and 15 isomorphically, hence these  
are also determined. In the 5th level starting at $F_{-1,4}=0$, the 
vanishing of the diagonal from top till $E^\infty_{4,1}=F_{4,1}/F_{3,2}$ 
forces the vanishing of the filtration till and including $F_{4,1}$. 
Now the limiting information $\mbbz_2=E^\infty_{5,0}=F_{5,0}/F_{4,1}$ 
determines the 5th homology. Similarly the 11th level can be handled. 
\end{proof}

\beg{lem}\label{H_4G} For the Grassmann manifold we have the following. 
$$H_5(G_3^+\mbb R^7;\mbbz)=H_6(G_3^+\mbb R^7;\mbbz)=\mbbz_2.$$
\end{lem}
\beg{proof} We saw that these two groups are solely torsion and isomorphic 
to one another and denoted both of them by $T_5$. We will be working on 
the cohomological Serre spectral sequence with integer coefficients 
related to the fiber
bundle \eqref{fibervg}, definition and limit of which is given as follows 
\cite{hajimesato}.
$$ E^{p,q}_2 := H^p( G \, ; H^q(SO_3;\mbb Z) ) $$
$$ E^{p,q}_\infty = { F^{p,q} / F^{p+1,q-1} }$$
where $F^{p,q}$ are abelian groups forming a filtration satisfying
$$H^n(V;\mbb Z) = F^{0,n} \supset F^{1,n-1} \supset \cdots \supset 
F^{n+1,-1}=0.$$
The differentials are of bidegree $(n,-n+1)$ so satisfying
$$d_n : E^{p,q}_n \longrightarrow E^{p+n,q-(n-1)}_n.$$

\beg{table}[ht]
\caption{ {\em Cohomological Serre spectral sequence for $G_3^+\mbb R^7$, 
second page.}}
\bct {\renewcommand*{\arraystretch}{1.2}
$\begin{array}{cc|c|c|c|c|c|c|c|} \cline{3-9}
&3&\,0\,& T_5 & T_5 &~~\,      &       &       &       \\ \cline{3-9}
&2&&    &T_{5,2}&\mbbz_2^2&\mbbz_2&       &       \\ \cline{3-9}
&1&&     &       &  \,0\, & \,0\, & \,0\, &       \\ \cline{3-9}
E_2~~~  &0&&&   &         &       &\mbbz_2& \,0\, \\ \cline{3-9}
& \multicolumn{1}{c}{ } 
& \multicolumn{1}{c}{5} & \multicolumn{1}{c}{6} 
& \multicolumn{1}{c}{7} & \multicolumn{1}{c}{8} 
& \multicolumn{1}{c}{9} & \multicolumn{1}{c}{10}
& \multicolumn{1}{c}{11}   \\  
\end{array}$  
} 
\ect
\setlength{\unitlength}{1mm}
\begin{picture}(0,0)
\put(77.5,30.5){\vector(3,-1){12}} 
\end{picture}
\label{table:E2Coh}
\end{table}

\noindent Our first claim is that $F^{7,3}\approx\mbbz_2$ for this sequence. Using the Theorem \ref{v37homology} and the filtration 
we obtain the isomorphisms
\beg{equation}\label{filtrationc}\mbbz_2=H_5\approx H^{10}(V_3\mbbr^7;\mbb Z) = F^{0,10} \supseteqq F^{1,9} \supseteqq \cdots \supseteqq F^{7,3}\end{equation}
provided by the vanishing of the limiting 
entries $E_\infty^{0,10}, \cdots, E_\infty^{6,4}$. Since we have 
$$\tn{Ker}d_2 = E_\infty^{7,3} = F^{7,3}/F^{8,2}$$
the only two possibilities $0$ or $\mbbz_2$ are remaining for $\tn{Ker}d_2$. After this point let us assume that \tcb{$T_5=0$} to raise a contradiction. 
Table \ref{table:E2CohHypothetical} shows the limit under this hypothesis.   
\beg{table}[ht]
\caption{ {\em 
Limiting page. Underlined terms are hypothetical.}}
\bct {\renewcommand*{\arraystretch}{1.2} 
$\begin{array}{cc|c|@{}c@{}|c|c|c|c|} \cline{3-8}
&3& ~~~ & \tn{Ker}d_2 &~~\,     &       &       &       \\ \cline{3-8}
&2&     & \tn{Ker}d_3 &\tcb{\underline{\mbbz}_2^2}& 0 && \\ \cline{3-8}
&1&     &     &         & \,0\, &       &       \\ \cline{3-8}
E_\infty~~~ &0&& & & &\tcb{\underline{\mbbz}_2}& ~~\\ \cline{3-8}
& \multicolumn{1}{c}{ } & \multicolumn{1}{c}{6} 
& \multicolumn{1}{c}{7} & \multicolumn{1}{c}{8} 
& \multicolumn{1}{c}{9} & \multicolumn{1}{c}{10}
& \multicolumn{1}{c}{11}   \\  
\end{array}$   
}
\ect
\label{table:E2CohHypothetical}
\end{table}
Note that the underlined terms are purely hypothetical. The two facts
$$\mbbz_2\approx F^{7,3} \supset F^{8,2} 
~~~\tn{and}~~~ 
\mbbz_2^2=E_\infty^{8,2}=F^{8,2}/F^{9,1}$$
shows that the group $\mbbz_2^2$ is way large to be carried by 
the filtration (\ref{filtrationc}). 
So that we now know $T_5\neq 0$. 
Next we claim that the torsion group $T_5$ is solely 2-torsion. 
To see it use the fundamental theorem of finitely generated abelian groups \cite{dummit} to conclude that this group is a direct sum of 
$\mbbz_{p^k}$'s for prime numbers $p\geq 2$ not necessarily distinct. If one of the $p$'s is odd than by the partial converse to Lagrange theorem (or the Sylow's theorem) there is a subgroup of order $p$. This subgroup is 
contained in the $\tn{Ker}d_2$ otherwise its image would be a group of order $2$ which has to divide $p$. However none of the two possibilities of $\tn{Ker}d_2$ above covers a subgroup of odd order. So $T_5$ is a direct sum of $\mbbz_{2^k}$'s. These summands are cyclic, so pick a generator i.e. an element so that the order $|a|=p^k$. 
Next we claim that $k$ cannot be greater than or equal to $3$. 
If that is the case to minimize the kernel $d_2$ must be surjective, 
in any case $|\tn{Ker}d_2|\geq 2^k/2=2^{k-1}\geq 2$ causes a problem. 
So that $k\leq 2$ hence $T_5$ can consist of $\mbbz_2$ or $\mbbz_4$ summands only. Computing the following entry of the spectral sequence 
 $$E_2^{7,2}=H^7(G;\mbbz_2)=\tn{Hom}(H_7,\mbbz_2)\oplus\tn{Ext}(H_6,\mbbz_2)=\tn{Ext}(T_5,\mbbz_2)$$ 
and denoting this term by $T_{5,2}$ (it counts the number of even ordered irreducible summands). We note that $T_{5,2}$ cannot hope to survive till infinity since $H^9(V;\mbbz)=F^{0,9}$ is trivial. So 
that the map 
$$d_3 : E^{7,2}_3=T_{5,2} \longrightarrow E^{10,0}_3=\mbbz_2$$
is injective. Reminding ourselves that $\tn{Ext}(\mbbz_4,\mbbz_2)=\tn{Ext}(\mbbz_2,\mbbz_2)=\mbbz_2$, the outcome is $T_{5,2}\subset\mbbz_2$. 
Since $T_5$ is nonzero we have $T_{5,2}=\mbbz_2$. 
Hence $T_5$ is either $\mbbz_2$ or $\mbbz_4$. 
To raise a contradiction, suppose $T_5=\mbbz_4$. 
Then $\tn{Ker}d_2$ cannot be zero by the cardinality, the remaining 
possibility is $\tn{Ker}d_2\approx \mbbz_2$ by above. 
Now we pass to the diagonal on the left. The differential 
$$d_2^{6,3} : \mbbz_4 \longrightarrow \mbbz_2\oplus\mbbz_2$$
cannot be surjective if it were, that would raise an isomorphism of the cyclic and Klein 4-group. So that it has a cokernel denoted 
${Cok}=E^{8,2}_\infty$ of order $2$ or $4$. 
Now concentrating on the $10th$ diagonal,  
our assumption  
$$\tn{Ker}d_2\approx \mbbz_2 = E^{7,3}_\infty=F^{7,3}/F^{8,2}$$
together with the fact that $F^{7,3}\approx\mbbz_2$ would imply 
that $F^{8,2}=0.$ Consequently we would obtain $E^{8,2}_\infty=0$, a contradiction.  
\end{proof}
\noindent We can now collect the results of this section to obtain the following. 
\beg{thm}\label{g37homology}The homology of the oriented Grassmann manifold $G_3^+\mbb R^7$ is computed as, 
$$H_*(G_3^+\mbb R^7;\mbb Z)=
(\mbbz,0,\mbbz_2,0,\mbbz\oplus\mbbz,\mbbz_2,\mbbz_2,
0,\mbbz\oplus\mbbz,\mbbz_2,0,0,\mbbz).$$
\end{thm}

\subsection{Cup product structure}\label{sectioncupproduct}

In this section we will analyze the cup product structure of the Grassmann 
manifold $G_3^+\mbbr^7$. We start with the free part. We will be using and interpreting the computations in \cite{zhou}. Recall $\xi$ and $\nu$ denote the 
{\em canonical (3-plane) bundle} and its orthogonal complement 4-plane bundle on this space respectively. Denoting the first Pontryagin and Euler classes of these bundles by $p=p_1(\xi)$ and $e=e(\nu)$ actually, we have the following. 

\beg{thm} The exterior algebra of the Grassmannian manifold $G_3^+\mbbr^7$ is given as follows $$H^*_{dR}(G_3^+\mbbr^7)=\mbbr[p,e]/(e^4,p^3-e^3,p^2-e^2,pe-ep) ~~\tn{for}~~ ~|p|=|e|=4.$$
\end{thm}
\beg{proof} Reading the Theorem 7.5 of \cite{zhou}, 
$p\pm e$ are generators of the fourth cohomology with Poincar\'e duals 
$2[ASS_{+}],2[ASS_{-}]$. And reading Theorem 7.4, $p^2, pe$ generate the eighth cohomology. The de Rham integral 
$\langle e^3, [G_3^+\mbbr^7]\rangle=2$ settles a non-zero class so that $e^3$ generates the top cohomology. So that the additive structure is given as follows. 
$$ H^*_{dR}(G_3^+\mbbr^7)=\< p+e,p-e\>\oplus \< p^2,pe \>
\oplus\< e^3 \>. $$
Among the relations $e^2=p^2$ is given in Section 7, 

$${ \renewcommand{\arraystretch}{2} 
\begin{array}{rcl} 
\< pe^2,[G_3^+\mbbr^7]\>  
& =& \< pe(e+p)/2 + pe(e-p)/2 , [G_3^+\mbbr^7] \>\\
& =& {1\over 2}\< pe\, , 2[ASS] \> + {1\over 2}\< pe\, , -2[\wt{ASS}] \>\\
& \stackrel{7.4}= & 2.\\
\end{array} }$$
\end{proof}

\beg{thm} The torsion algebra of the Grassmannian manifold $G_3^+\mbbr^7$ is given as follows 
$$\mbbz_2[x_3,x_7]/(x_3^3,x_7^2,x_3x_7+x_7x_3) ~~\tn{for}~~ ~|x_3|=3, |x_7|=7.$$
\end{thm}
\beg{proof} There is torsion at four levels $3,6,7,10$ as we have computed in section \ref{sectiongrassmanninvariants}.  
To understand the cup product structure we consult to the cohomological Serre spectral sequence with $\mbbz$-coefficients. 
\beg{table}[ht]
\caption{\em Cohomological Serre spectral sequence for $G_3^+\mbb R^7$. 
Second page.} \bct{
{\renewcommand*{\arraystretch}{1.2}    
$\begin{array}{cc|c|c|c|c|c|c|c|c|c|c|c|c|c|} \cline{3-15}
&3&\mbbz&\,0\,&\,0\,&\mbbz_2&\mbbz^2&\,0\,&\mbbz_2&\mbbz_2&\mbbz^2&\,0\,&\mbbz_2&\,0\,&\mbbz \\ \cline{3-15}
&2&\mbbz_2&\,0\,&\mbbz_2&\mbbz_2&\mbbz_2^2&\mbbz_2&\mbbz_2^2&\mbbz_2&
\mbbz_2^2&\mbbz_2&\mbbz_2&\,0\,&\mbbz_2 \\ \cline{3-15}
&1&\,0\,&\,0\,&\,0\,&\,0\,&\,0\,&\,0\,&\,0\,&\,0\,&\,0\,&\,0\,&\,0\,&\,0\,&\,0\,      \\ \cline{3-15}   E_2~~~
&0&\mbbz&\,0\,&\,0\,&\mbbz_2&\mbbz^2&\,0\,&\mbbz_2 &\mbbz_2 &\mbbz^2&0&\mbbz_2&\,0\,&\mbbz \\ \cline{3-15}
& \multicolumn{1}{c}{ } & \multicolumn{1}{c}{0} 
& \multicolumn{1}{c}{1} & \multicolumn{1}{c}{2} 
& \multicolumn{1}{c}{3} & \multicolumn{1}{c}{4}
& \multicolumn{1}{c}{5} & \multicolumn{1}{c}{6} 
& \multicolumn{1}{c}{7} & \multicolumn{1}{c}{8} 
& \multicolumn{1}{c}{9} & \multicolumn{1}{c}{10}   
& \multicolumn{1}{c}{11} & \multicolumn{1}{c}{12}   
\\  
\end{array}$   
}   } \ect
\label{table:E2cohom}
\end{table}
Table \ref{table:E2cohom} shows the second third page of this spectral
sequence. 
In Table \ref{table:E3cohom} for the third page of the spectral
sequence, 
the arrows 
are isomorphism as follows.
We have $H^3(V;\mbbz)=0=F^{0,3}$ from Theorem \ref{v37homology} implying that the limit term $E_\infty^{3,0}=F^{3,0}/F^{4,-1}$ vanishes. To provide that the only incoming 
non-zero differential $d_3:E_3^{0,2}\to E_3^{3,0}$ must be an isomorphism. 
The second one is similar, the vanishing of $E_\infty^{3,2}=F^{3,2}/F^{4,1}$ is guaranteed through the vanishing of 
$H^5(V;\mbbz)=F^{0,5}$. The last one is achieved through $H^9(V;\mbbz)=0=F^{0,9}$ hence $E_\infty^{7,2}=0$. 
\beg{table}[ht]
\caption{\em Cohomological Serre spectral sequence for $G_3^+\mbb R^7$. 
} \bct{ 
{\renewcommand*{\arraystretch}{1.2}
$\begin{array}{cc|c|c|c|@{}c@{}|c|c|@{}c@{}| @{}c@{}|c|c| @{}c@{} |} \cline{3-13}
&2&\mbbz_2a&~~&~~&\mbbz_2ax_3&\mbbz_2^2&~~&0&\mbbz_2ax_7&~~&~~& \\ \cline{3-13}
&1&       &&&&&         &&&&&      \\ \cline{3-13}   E_3~~~
&0&\mbbz 1&&&\mbbz_2x_3&&&\mbbz_2 x_6&\mbbz_2x_7&&0&\mbbz_2 x_{10}\\ \cline{3-13}
& \multicolumn{1}{c}{ } & \multicolumn{1}{c}{0} 
& \multicolumn{1}{c}{1} & \multicolumn{1}{c}{2} 
& \multicolumn{1}{c}{3} & \multicolumn{1}{c}{4}
& \multicolumn{1}{c}{5} & \multicolumn{1}{c}{6} 
& \multicolumn{1}{c}{7} & \multicolumn{1}{c}{8} 
& \multicolumn{1}{c}{9} & \multicolumn{1}{c}{10}   
\\  
\end{array}$   
\setlength{\unitlength}{1mm}
\begin{picture}(0,0)
\put(-84,7){\vector(2,-1){13}}   
\put(-60,7){\vector(2,-1){13}}
\put(-26,7){\vector(2,-1){13}}   
\end{picture} }
}
\ect
\label{table:E3cohom}
\end{table}

\ni Following the techniques in \cite{hatcherspec} we label the generators as above. Replacing some generators with their negatives if necessary we may assume $d_3a=x_3$. Similarly we may assume $d_3(ax_3)=x_6$. Combining with the following relation  
$$d_3(a x_3)=d_3ax_3+adx_3=x_3^2+a0=x_3^2$$
we replace $x_6=x_3^2$. Again assuming $d_3(ax_7)=x_{10}$ and applying the following identity 
$$d_3(a x_7)=d_3ax_7+adx_7=x_3x_7$$ 
we replace $x_{10} = x_3 x_7$. Since there is no cohomology at the level $9$, $x_3^3=0$. One can alternatively see this as follows. 
Since $H^8(V;\mbbz)=0=F^{0,8}$ we have $E^{6,2}_\infty=0$ forces that the entry $E^{6,2}_3=0$. From the isomorphism 
$E_3^{0,2}\otimes E_3^{6,0}\longrightarrow E_3^{6,2}$ we get 
$ax_6=ax_3^2=0$. Taking the differential of both sides yields the result. 
$$0=d_3(ax_3^2)=d_3(a)x_3^2+ad_3(x_3^2)=x_3^3.$$ 
Last relation follows from the alternating property of the cup product for odd dimensions.
\end{proof}

\ni Combining the two results we obtain the following. 
\beg{thm}\label{cupproductstructuregrassmannian} The cohomology ring of the Grassmannian manifold $G_3^+\mbbr^7$ is given as follows 
$$H^*(G_3^+\mbbr^7;\mbbz)=\mbbz[x,y]/(x^4,x^3-y^3,x^2-y^2,xy-yx) 
\oplus\mbbz_2[z,t]/(z^3,t^2,zt+tz)$$ 
where the degrees are~ $|x|=|y|=4, |z|=3$ ~and~ $|t|=7$, and $x=e(\nu)$, $y=p_{1}(\xi)$.
\end{thm}

Next we would like to see how the submanifold of associative planes sits 
inside $G_3^+\mbb R^7$ cohomologically. We would like to mention that the space of associative 3-
planes and its Stiefel-Whitney classes are studied to some degree by \cite{borelhirzebruch1}. We 
have the pullback map induced by the inclusion
$$H^q(ASS;\mbbz) \longleftarrow H^q(G_3^+\mbb R^7;\mbbz) : i^*$$
which operates at the levels $q=0 \cdot\cdot\, 8$. The nonzero integral cohomology groups are 
already computed in Section 10 of \cite{zhou} to be the following. 
$$H^q(ASS;\mbbz)=\left\{
\beg{array}{ll}
\mbbz,  & q=0,4,8\\
\mbbz_2,& q=3,6\\
\end{array}
\right.$$
We can compute the cohomological ring and the action of the inclusion map on cohomology as follows. 

\beg{thm} We have the following facts for the $8$-manifold of associative planes.
\beg{enumerate}
\item[(a)] The cohomology ring structure is given as
$$H^*(ASS;\mbbz)=\mbbz[d]/\<d^3\>\oplus\mbbz_2[c]/\<c^3\>$$ 
where the degrees are~ $|d|=4$\, and\, $|c|=3$.
\item[(b)] The inclusion map $i: ASS \to G_3^+\mbbr^7$ acts on the cohomology rings as follows
$$i^*x=i^*y=d, ~~i^*z=c,~~i^*t=0.$$
\end{enumerate}
\end{thm}
\begin{proof} The proof that we will give for the two parts are somehow interrelated. \begin{enumerate}
\item We will employ the fibration $S^3\to G_2^+\mbbr^7\to ASS$, 
see \cite{zhou} for the map. Corresponding Serre cohomological spectral sequence yields the following page. 
\beg{table}[ht]
\caption{\em Cohomological Serre spectral sequence for $G_2
^+\mbb R^7$. 
} \bct{ 
{\renewcommand*{\arraystretch}{1.2}
$\begin{array}{cc|c|c|c|@{}c@{}|@{}c@{}|c|@{}c@{}|c|c|} \cline{3-11}
&2&\mbbz a&~~~&~~~&\mbbz_2ax_3&\mbbz ax_4&~~~&\mbbz_2ax_6&~~~&\mbbz ax_8 \\ \cline{3-11}
&1&       &&&&&         &&&      \\ \cline{3-11}   E_3=E_2~~~
&0&\mbbz 1&&&\mbbz_2x_3&\mbbz x_4&&\mbbz_2 x_6&&\mbbz x_8 \\ \cline{3-11}
& \multicolumn{1}{c}{ } & \multicolumn{1}{c}{0} 
& \multicolumn{1}{c}{1} & \multicolumn{1}{c}{2} 
& \multicolumn{1}{c}{3} & \multicolumn{1}{c}{4}
& \multicolumn{1}{c}{5} & \multicolumn{1}{c}{6} 
& \multicolumn{1}{c}{7} & \multicolumn{1}{c}{8} 
\\  
\end{array}$   
\setlength{\unitlength}{1mm}
\begin{picture}(0,0)
\put(-75,7){\vector(2,-1){13}}   
\put(-49,7){\vector(2,-1){13}}
\put(-37,7){\vector(2,-1){13}}   
\end{picture} }
}
\ect
\label{table:associatesimplestfibration}
\end{table}
\ni Here $c$ appears as the pullback of the Euler class of the tautological bundle on the Grassmannian which is non-zero. Labeling 
the generators in the spectral sequence as in Table \ref{table:associatesimplestfibration}, the surjective map 
yields the relation $d_3a=x_3$. Then the isomorphism implies 
$x_6=d_3(ax_3)=d_3a x_3\pm adx_3=x_3^2$. This shows that $c^2$ 
is also a nonzero element, hence the generator of its level. 
Another relation is obtained through the third map $0=d_3(ax_4)=x_3x_4$. 
Relabel $d=x_4$.
\item Next we will deduce that the element $d^2$ is the generator of its 
level. And this will finish the proof on the part $(a)$. 
To see this let 
$$H^4(ASS;\mbbz)=\<d\> ~~\tn{and}~~ H^8(ASS;\mbbz)=\<s\>.$$ 
Then $d^2=\alpha\, s$ for some $\alpha\in\mbbz$. Also suppose that 
$i^*x=\beta\, d$ and $i^*y=\gamma\, d$ for some $\beta,\gamma\in\mbbz$. 
Recall that additively we have the following generators at 
the level $8$ of the Grassmannian, 
$$H^8(G_3^+\mbbr^7;\mbbz)=\<p_1^2E,p_1EeF\>$$
where $x=p_1E$ and $y=eF$. By the Lemma 7.4 of \cite{zhou} we have 
the following integrals 
$$\<x^2,[ASS]\>=\<xy,[ASS]\>=1,$$
so that these elements map onto the generator $d^2$ of the eighth 
cohomology of $ASS$. Elaborating 
this fact by {\renewcommand*{\arraystretch}{1.2}
$$\begin{array}{cccccc} 
s&=& i^*x^2     & = & i^*(xy)    \\
 &=& \beta^2\,d^2 & = & \gamma^2\,d^2
\end{array}$$}
\ni and plugging in yields 
$$s= \beta^2 \alpha \,s.$$
So that $\beta^2 \alpha=1$. This forces\, $\alpha=1$ hence the assertion.

\item Since now we have obtained $s=d^2$, we have $|\beta|=|\gamma|=1$. 
We can assume that $\beta=1$ after a change of 
sign of the generator $d$ if needed. We claim the same for $\gamma$ 
as well, suppose $\gamma=-1$ to raise a contradiction. Then $i^*(x+y)=0$ 
and this implies 
 {\renewcommand*{\arraystretch}{1.2}
$$\begin{array}{ccl} 
0&=& i^*(x+y)^2     \\
 &=& i^*(x^2+xy+yx+y^2)\\
 &=& d^2+2i^*(xy)+d^2 
\end{array}$$}
\ni yields the contradiction $i^*(xy)=-d^2$ to facts of the previous part. 
Combining these we have 
$i^*x=i^*y=d$.
\item The element $c$ comes naturally as the restriction of the Euler 
class $z$ of the tautological vector bundle $E_3^7$. Finally there is no 
seventh cohomology to map onto.\end{enumerate}\end{proof}


\subsection{The Lie Group $G_2$}\label{sectiong2homology}
In this section we will compute some invariants of the Lie Group $G_2$, which is defined to be the subgroup of $SO_7$ which fixes the 4-form,
$$*\phi_0=dx^{4567}+dx^{2367}+dx^{2345}+dx^{1357}-dx^{1346}-dx^{1256}-dx^{1247}.$$
Here the notation $dx^{4567}$ suggest $dx^4\wedge dx^5\wedge dx^6\wedge dx^7$ similarly the 
others. In order to work efficiently on $G_2$ we will need to use two fibrations first of which 
is the following. 
\beg{equation}\label{g2su3fibration}SU_3\longrightarrow G_2 \longrightarrow S^6\end{equation}
Since $G_2$ consists of orthogonal transformations it preserves the 
sphere in $\mbb R^7$ so that it has an action on the 6-sphere.  
This is a transitive action and the stabilizer of a point on the sphere,
preserves its orthogonal complement as well hence a subgroup of $SO_6$. 
See a general reference \cite{bryant} for further details. 

\vspace{.05in}

To work with the fibration we need the cohomology of the fiber. 
One can obtain the cohomology of the unitary group as an exterior algebra 
$$H^*(U_n;\mbbz)=\Lambda_\mbbz[x_1,x_3\cdots x_{2n-1}]$$
using the complex Stiefel manifolds. Then via the action of the special unitary group on the unitary one, the fibration 
$SU_n\to U_n\to S^1$ helps to drop the first generator and we get 
$$H^*(SU_3;\mbbz)=\Lambda_\mbbz[x_3,x_5]=(\mbbz,0,0,\mbbz,0,\mbbz,0,0,\mbbz).$$
Now we can work on the cohomological Serre spectral sequence for the first fibration (\ref{g2su3fibration}) which up to the sixth page looks like in Table \ref{table:g2su3cohom}.
\beg{table}[ht]
\caption{ {\em Cohomological Serre spectral sequence for the $SU_3$ fibration of $G_2$.
}} \bct{
$\begin{array}{cc|c|c|c|c|c|c|c|} \cline{3-9}
  &8& \mbbz &&&&&& \mbbz\\ \cline{3-9}
  &7&       &&&&&&      \\ \cline{3-9}
  &6&       &&&&&&      \\ \cline{3-9}
  &5& \mbbz &&&&&& \mbbz\\ \cline{3-9}
  &4&       &&&&&&      \\ \cline{3-9}
  &3& \mbbz &&&&&& \mbbz\\ \cline{3-9}
  &2&       &&&&&&      \\ \cline{3-9}
  &1&       &&&&&&      \\ \cline{3-9}          E_2=E_6~~~
  &0& \mbbz &&&&&& \mbbz\\ \cline{3-9}
& \multicolumn{1}{c}{ } & \multicolumn{1}{c}{0} 
& \multicolumn{1}{c}{1} & \multicolumn{1}{c}{2} 
& \multicolumn{1}{c}{3} & \multicolumn{1}{c}{4}
& \multicolumn{1}{c}{5} & \multicolumn{1}{c}{6} 
\\  
\end{array}$   
\setlength{\unitlength}{1mm}
\begin{picture}(0,0)
\put(-38,24){\vector(4,-3){30}}   
\put(-38,8 ){\vector(4,-3){30}}
\end{picture}
}\ect
\label{table:g2su3cohom}
\end{table}
Note that the real dimension of the Lie algebra $\mathfrak{su}_3$ is
computed to be 8 so that $G_2$ is 14 dimensional.   We can immediately compute from this sequence that there is no cohomology at the levels 
$1,2,4,7,10,12,13$  and there is a $\mbbz$ each at the levels $3,11$. 
This is more or less the only accessible information to get at first sight from this fibration, this is mainly because we do not use the 
actual fibration map which could have been the trivial product $SU_3\times S^6$ as well. Since this is a fibration the Euler characteristic is multiplicative and 
$$\chi(G_2)=\chi(SU_3)\chi(S^6)=0\cdot 2$$  
so that this implies $b_8=b_5$. At this point using universal coefficients and Poincar\'e duality,$$H^5=H^8=F_5~~\tn{and}~~H^6=H^9=F_5\oplus T_5.$$
is the ultimate statement for the missing cohomology along with the following Lemma. 
Here $F_5$ and $T_5$ denotes the free and torsion part of the fifth homology.   


\beg{lem}\label{b5leq1} The fifth Betti number $b_5(G_2)\leq 1$.
\end{lem}
\beg{proof} From the filtration one can show that 
$$H^9(G_2;\mbbz)\approx F^{6,3}\approx E^{6,3}_\infty=\mbbz/\tn{Im}\,d^{0,8}_6.$$
A similar computation can be done to see $H^6(G_2;\mbbz)\approx\mbbz/\tn{Im}\,d^{0,5}_6.$ 
Since we know that these two are isomorphic, as a by product one can easily see the equality $\tn{Im}\,d^{0,8}_6=\tn{Im}\,d^{0,5}_6$. 
\end{proof}

To recover the missing information we consult to the second  fibration (\ref{g2so4fibration}). To analyze this we need the homology of the fiber. We use the fibration 
\beg{equation}\label{so4so3fibration}SO_3\longrightarrow SO_4 \longrightarrow S^3\end{equation}
and applying the homological spectral sequence as in the Table \ref{table:so4so3hom} yields the homology,
$$H_*(SO_4;\mbbz)=(\mbbz,\mbbz_2,0,\mbbz^2,\mbbz_2,0,\mbbz).$$
Here one should use universal coefficients to compute the homology at level 4 and to show it is free at level 3. Then computing the Euler characteristic $\chi(SO_4)=\chi(SO_3)\chi(S^3)$ to be zero yields $b_3=2$,
hence the result. 

{ \renewcommand*{\arraystretch}{1.1} 
\beg{table}[ht]
\caption{ {\em Homological Serre spectral sequence for the $SO_3$ fibration of $SO_4$. 
}}
\bct{
$\begin{array}{cc|c|c|c|c|c|} \cline{3-6}
  &3& \mbbz & \mbbz_2 &\,0\,& \mbbz \\ \cline{3-6}
  &2&   0   &    0    &  0  &   0   \\ \cline{3-6}
  &1&\mbbz_2&    0    &  0  &\mbbz_2\\ \cline{3-6}          E^2=E^\infty~~~
  &0& \mbbz & \mbbz_2 &  0  & \mbbz \\ \cline{3-6}
& \multicolumn{1}{c}{ } & \multicolumn{1}{c}{0} 
& \multicolumn{1}{c}{1} & \multicolumn{1}{c}{2} 
& \multicolumn{1}{c}{3}
\\  
\end{array}$   
\setlength{\unitlength}{1mm}
\begin{picture}(20,10)
\put(-16,-3){\vector(-3,1){12}}
\end{picture}
}\ect
\label{table:so4so3hom}
\end{table}  }

\ni The second fibration of $G_2$ is as follows.
\beg{equation}\label{g2so4fibration}SO_4\longrightarrow G_2 \longrightarrow ASS\end{equation}
To see this, note that the group $G_2\subset SO_7$ naturally acts on $\mbb R^7$, leaves the associative form invariant. So that an associative (three) plane is sent to another associative plane under $G_2$ action, hence the set ASS of associative planes stays invariant. The stabilizer 
of this action is the orthogonal transformations of $\mbb R^7$ which 
leaves an associative 3-plane invariant, and hence acts orthogonally in 
the complement yielding $SO_4$. 
{\renewcommand*{\arraystretch}{1.2}    
\beg{table}[ht]
\caption{ {\em Homological Serre spectral sequence for the $SO_4$ fibration of $G_2$.
}}
\bct{
$\begin{array}{cc|c|c|c|c|c|c|c|c|c|} \cline{3-11}
&6& \mbbz &0&\mbbz_2&0&\mbbz&\mbbz_2&0&0&\mbbz_2          \\ \cline{3-11}
&5&   0   &\,0\,&\,0\,&\,0\,&\,0\,&\,0\,&\,0\,&\,0\,&\,0\,\\ \cline{3-11}
&4&\mbbz_2&0&\mbbz_2&\mbbz_2&\mbbz_2&\mbbz_2&\mbbz_2&0&\mbbz_2\\\cline{3-11}
&3&\mbbz^2&0&\mbbz_2^2&0&\mbbz^2&\mbbz_2^2&0&0&\mbbz^2    \\ \cline{3-11}
&2&   0   &\,0\,&\,0\,&\,0\,&\,0\,&\,0\,&\,0\,&\,0\,&\,0\,\\ \cline{3-11}
&1&\mbbz_2&0&\mbbz_2&\mbbz_2&\mbbz_2&\mbbz_2&\mbbz_2&0&\mbbz_2\\\cline{3-11}          E^2~~~
&0& \mbbz &0&\mbbz_2&0&\mbbz&\mbbz_2&0&0&\mbbz        \\ \cline{3-11}
& \multicolumn{1}{c}{ } & \multicolumn{1}{c}{0} 
& \multicolumn{1}{c}{1} & \multicolumn{1}{c}{2} 
& \multicolumn{1}{c}{3} & \multicolumn{1}{c}{4}
& \multicolumn{1}{c}{5} & \multicolumn{1}{c}{6} 
& \multicolumn{1}{c}{7} & \multicolumn{1}{c}{8} 
\\  
\end{array}$   
\setlength{\unitlength}{1mm}
\begin{picture}(0,0)
\put(-60,-14){\vector(-2,1){10}}   
\end{picture}
}\ect
\label{table:g2so4hom}
\end{table}   }
According to the rule 
$$E^2_{p,q}=H_p(ASS;H_q(SO_4;\mbbz))$$
Table \ref{table:g2so4hom} shows the starting page for the homological Serre spectral sequence for this fibration. The following is our first main assertion. 
\beg{lem} The free part $F_5$ is zero hence $b_5(G_2)=0$.
\end{lem}
\beg{proof} We actually claim that the line $p+q=9$ at infinity is totally zero. To see this consider the filtration
$$0=F_{-1,10} =
    F_{ 0, 9} =
\cdots
   =F_{ 4, 5} \subset  
\cdots        \subset 
    F_{7,2}   \subset 
    F_{8,1}   =
    F_{9,0}=H_9(G_2;\mbbz)=F_5.$$
Over the line at infinity the only possibly nonzero terms are 
$E_{5,4}^\infty$ and $E_{8,1}^\infty$ which both are subgroups of $\mbbz_2$. Since lower terms at infinity are zero we have 
$E_{5,4}^\infty=F_{5,4}$. So that $F_{5,4}$ is a subgroup of a free group and a subgroup of $\mbbz_2$ at the same time. So it has to vanish. 
Then if you follow up zeros till $(8,1)$ you can do the same argument to see that $F_{8,1}=0$, hence the result.  \end{proof}
This gives us a chance to say something about the torsion. 
\beg{lem} The torsion part $T_5=\mbbz_m$ for some $m\geq 2$ for the Lie group $G_2$. In particular it is nonzero. 
\end{lem}
\beg{proof} Recall from Lemma \ref{b5leq1} that we have
$$H^9(G_2;\mbbz)\approx \mbbz/\tn{Im}\,d^{0,8}_6=F_5+T_5=T_5.$$
Assume that $d^{0,8}_6$ is surjective or $T_5=0$ to raise a contradiction. 
In that case 
passing to the second (homological) sequence where 
$H_8(G_2;\mbbz)=0$ and the $p+q=8$ line disappears in the limit. 
In particular $E_{5,3}^\infty=0$. 
The only differential from or hitting $(5,3)$ are  
$d_{5,3}^2$ and $d_{8,1}^3$. 
$d_{5,3}^2$ should better be surjective 
and $d_{8,1}^3$ better be an embedding to bleed $E_{5,3}^2=\mbbz_2^2$ to nothing since these two are the only two chances.   
The second assertion means $$E_{8,1}^4:=\tn{Ker}\,d_{8,1}^3=0.$$ 
But now consider the unpleasant situation for $E_{4,4}^2=\mbbz_2$ which has to 
bleed into death. Only possibly nontrivial differential is the following 
$$d_{8,1}^4: E_{8,1}^4 \longrightarrow E_{4,4}^4$$
which emanates from zero as we computed, a contradiction. \end{proof}

We will also be using the following Lemma. 

\beg{lem}\label{tepeisomorfizm} 
We have that $E_{2,6}^\infty=E_{8,0}^\infty=0$ for the limits. 
Moreover $E_{5,4}^3=E_{2,6}^3=0$. 
Hence only possibly nonzero terms on the line $p+q=8$ at infinity are 
$E_{4,4}^\infty$ and $E_{5,3}^\infty$.
\end{lem}
\beg{proof} 
Since we know that $E_{5,4}^\infty=0$,  the differential $$d_{5,4}^2 : E_{5,4}^2 \longrightarrow E_{2,6}^2$$
has to be injective hence an isomorphism, enough to kill the entry $(2,6)$. Moreover this implies also that $E_{5,4}^3=E_{2,6}^3=0$ on both parts. 
Since $E_{8,0}^2=\mbbz$ is free, so is $E_{8,0}^\infty=F_{8,0}/F_{7,1}$. 
Considering $F_{8,0}=H_8(G_2;\mbbz)=T_5=\mbbz_m$,  $E_{8,0}^\infty$ also torsion hence trivial. \end{proof}

Finally we can now handle the torsion piece.

\beg{lem} The torsion part $T_5=\mbbz_2$ for the Lie group $G_2$.
\end{lem}
\beg{proof} 
All of the terms from  $E_{0,8}^\infty$ till and including  $E_{3,5}^\infty$ vanish. On the filtration this implies that 
$$0=F_{-1,9} =
    F_{ 0, 8} =
    F_{ 1, 7} =
    F_{ 2, 6} =
    F_{ 3, 5} \subset  
    F_{ 4, 4} \subset  
\cdots        \subset 
    F_{8,0}=H_8(G_2;\mbbz)=T_5=\mbbz_m$$
the left hand  terms vanish and $E_{4,4}^\infty\approx F_{4,4}$. From the starting entry of the spectral sequence we have 
$E_{4,4}^\infty\unlhd \mbbz_2$ hence $F_{4,4}\unlhd \mbbz_2$ as well. 
So we have only two possibilities for $F_{4,4}$. 
We will analyze these two cases separately.

\beg{itemize}

\item[{Case 1:}] Assume $F_{4,4}=\mbbz_2$.  
Then all differentials related to the $(4,4)$ terms are zero for its survival. In particular the differential 
$$d_{8,1}^4: E_{8,1}^4 \longrightarrow E_{4,4}^4$$
is zero. That implies the convergence $E_{8,1}^4=E_{8,1}^\infty=0$. 
Since from Lemma \ref{tepeisomorfizm} we have $0=E_{2,6}^3=E_{2,6}^6$ so that the differential 
$d_{8,1}^6:E_{8,1}^6\longrightarrow E_{2,6}^6$ is zero. Consequently the only nonzero differential emanating from the entry $(8,1)$ is the 
\beg{equation}\label{differential}
d_{8,1}^3:E_{8,1}^3=\mbbz_2\longrightarrow E_{5,3}^3=K
\end{equation}
which has to be injective to provide $E_{8,1}^\infty=0$. Here, by definition we take $K:=\tn{Ker}\,d_{5,3}^2$ so that this kernel has a subgroup of order two, in particular it is nontrivial. In the following we
 will show that on the other end $K\neq\mbbz_2^2$ as well. To see this 
 observe that the differential 
$$d_{8,0}^5:E_{8,0}^5\lhd\mbbz\longrightarrow E_{3,4}^5\lhd \mbbz_2$$
has to be zero. If it does not then it would mean that $E_{8,0}^5\approx\mbbz$. But this term vanishes ultimately 
by Lemma \ref{tepeisomorfizm} and there 
is no chance to vanish since the differential can no longer be injective.
The outgoing differential $d_{3,4}^3 : \mbbz_2 \longrightarrow \mbbz$ is zero. So the only possibly nontrivial differential concerning $(3,4)$ 
is $$d_{5,3}^2:E_{5,3}^5=\mbbz_2^2\longrightarrow E_{3,4}^2=\mbbz_2$$
has to be surjective to kill it since $H_7(G_2;\mbbz)=0$ so that it disappears at infinity. To be surjective the kernel $K$ cannot be everything. 
So we reach at the only possibility that $K\approx\mbbz_2$. 
But the differential at (\ref{differential}) was injective 
now becomes an isomorphism. That kills the term and we get $E_{5,3}^\infty=0$. That tells the result as follows
$$\mbbz_2=F_{4,4}=F_{5,3}=F_{6,2}=F_{7,1}=F_{8,0}=\mbbz_m.$$

\item[{Case 2:}]  Assume $F_{4,4}=0$. Then $E_{4,4}^\infty=F_{4,4}/F_{3,5}=0$. 
Then the only possibly nontrivial differential incoming or emanating from 
$E_{4,4}^*$ is 
$$d_{8,1}^4: E_{8,1}^4 \longrightarrow E_{4,4}^4=\mbbz_2$$
which has to be surjective. So $E_{8,1}^4=\mbbz_2$, implying that 
$E_{8,1}^3=\mbbz_2$ as well and the map 
$$d_{8,1}^4: E_{8,1}^3=\mbbz_2 \longrightarrow E_{5,3}^3$$
is zero. Letting $K:=\tn{Ker}\,d_{5,3}^2$ then we have 
$E_{5,3}^\infty=K$. The elements on the filtration become 
$$F_{5,3}=F_{6,2}=F_{7,1}=F_{8,0}=\mbbz_m,$$
so that $$K=E_{5,3}^\infty=F_{5,3}/F_{4,4}=F_{5,3}=\mbbz_m.$$
We also know that $K\unlhd \mbbz_2^2$ so that $\mbbz_m \unlhd \mbbz_2^2$ and 
it $m\geq 2$ implies that $m=2$.
\end{itemize}
\end{proof}

\ni Accumulating the results of this section we get the following. 
\begin{thm}\label{g2homology} 
The homology groups of the Lie group $G_2$ are as follows.
$$H_*(G_2;\mbb Z)=
(\mbbz,0,0,\mbbz,0,\mbbz_2,0,0,\mbbz_2,0,0,\mbbz,0,0,\mbbz).$$
\end{thm}

At this point we are in a position to understand the cup product structure of $G_2$ as follows.

\begin{thm}\label{g2cup} 
The cohomology ring of the Lie group $G_2$ can be described as follows.
$$H^*(G_2;\mbb Z)=\Lambda_\mbbz[x_3,x_{11}]\oplus
\Lambda_{\mbbz_2}[x_6,x_9]/(x_6 x_9)$$
where the degrees are~ $|x_k|=k$.
\end{thm}

\beg{proof} Consider Table \ref{table:g2su3cup} that presents all possible generators. 
{\renewcommand*{\arraystretch}{1.2} 
\beg{table}[ht]
\caption{ {\em Cohomological Serre spectral sequence for the $SU_3$ fibration of $G_2$.
}} \bct{
$\begin{array}{cc|c|ccccc|c|} \cline{3-9}
&8&\mbbz x_8&&&&&& \mbbz x_{14} \\ \cline{3-9}
&7&         &&&\times 2&&&              \\ 
&6&         &&&&&&              \\ \cline{3-9}
&5&\mbbz x_5&&&&&& \mbbz x_{11} \\ \cline{3-9}
&4&         &&&\times 2&&&              \\ \cline{3-9}
&3&\mbbz x_3&&&&&& \mbbz x_9    \\ \cline{3-9}
&2&         &&&&&&              \\ 
&1&         &&&&&&              \\ \cline{3-9} E_2=E_6~~~
&0&\mbbz 1  &&&&&& \mbbz x_6    \\ \cline{3-9}
& \multicolumn{1}{c}{ } & \multicolumn{1}{c}{0} 
& \multicolumn{1}{c}{1} & \multicolumn{1}{c}{2} 
& \multicolumn{1}{c}{3} & \multicolumn{1}{c}{4}
& \multicolumn{1}{c}{5} & \multicolumn{1}{c}{6} 
\\  
\end{array}$   
\setlength{\unitlength}{1mm}
\begin{picture}(0,0)
\put(-45,24){\vector(4,-3){30}}   
\put(-45,6 ){\vector(4,-3){30}}
\end{picture}
}\ect \label{table:g2su3cup}
\end{table} }
The maps are multiplication by $2$ which are injective and so that 
$x_5,x_8$ disappears.  
Since the image is zero, the product
$$E^{0,3}_6\otimes E^{0,3}_6 \longrightarrow E^{0,6}_6 $$
is zero which is the same as cup product upto a sign. So that this implies 
$x_3^2=0$. Likewise one can compute $x_6^2=x_9^2=x_{11}^2=0$. 
Since $$E^{0,3}_6=H^0(S^6;H^3(SU_3;\mbbz))\approx H^3(SU_3;\mbbz)$$
the product
$$E^{0,3}_6\otimes E^{6,5}_6 \longrightarrow E^{6,8}_6 $$
is just multiplication of coefficients. Hence the multiplication 
$x_3\cdot : E^{6,5}_6 \longrightarrow E^{6,8}_6 $ is an isomorphism. 
In particular sends generators to generators hence $x_3\cdot x_{11}=x_{14}$ if the signs are arranged suitably. The only missing relation 
$$x_3x_{11}=(-1)^{3\cdot 11} x_{11}x_3$$ 
comes by the properties of the cup product.
\end{proof}

\subsection{Classifying Space}\label{sectionclassifying}
In this section we will compute the classifying space $BG_2$ of the group 
$G_2$. We will be using the results \cite{borel} of A. Borel. See 
\cite{mimura} for a recent exposition. We start with the free algebra. 
A theorem of Borel tells us the following. Let $G$ be a compact, connected Lie group and $R$ be $\mbbz$ or a field $k$ of characteristic $p$. 
Assume $H_*(G)$ is torsion free if $R=\mbbz$, or is $p$-torsion free 
if $R=k$. Then there are universally transgressive elements 
$x_i\in H^{n_i}(G;R)$ such that 
$$H^*(G;R)=\Lambda[x_1\cdots x_l]~~\tn{where}~~ |x_i|=n_i~\tn{odd},$$   
$$H^*(BG;R)=R[y_1\cdots y_l]~~\tn{where}~~ y_i=\tau(x_i).~~~~$$

\ni Here $\tau: H^k(BG)\to H^{k+1}(G)$ is the transgression map of the fiber bundle $G\to EG\to BG$. This is a co-analogue of the connecting homomorphism of the homotopy exact sequence of a fiber bundle. A corollary of this theorem is that if $H_*(G)$ is $p$-torsion free, then $H_*(BG)$ 
is also $p$-torsion free. In our case this means that $BG_2$ has no $p$-torsion for $p\geq 3$. At this point, taking a $\mbbz_3$ coefficient 
cohomology ring kills the ($2$) torsion and captures the free piece by universal coefficients theorem. Starting with the free part of $G_2$ 
as we computed in Theorem \ref{g2cup} we have generators at the levels $3$ 
and $11$. Transgression increases the degree by one and applying the above we have 
$$H^*(BG_2;\mbbz_3)=\mbbz_3[y_4,y_{12}].$$
For the torsion part, by the application of the above theorems, 
it is well-known that the cohomology with $\mbbz_2$ coefficients is 
$$H^*(G_2;\mbbz_2)=\mbbz_2[x_3]/(x_3^4)\otimes \Lambda[x_5]~~\tn{where}~~
x_5=Sq^2x_3,$$
$$H^*(BG_2;\mbbz_2)=\mbbz_2[y_4,y_6,y_{10}]~~\tn{where}~~
y_6=Sq^2y_4~~\tn{and}~~y_7=Sq^3y_4=Sq^1y_6.$$
Now, comparing the generators, in this $\mbbz_2$ coefficient ring, 
$y_4$ comes from the free part, $y_6,y_{10}$ are new so that they are produced by the torsion. 
\beg{thm}\label{bg2}
The cohomology ring of the classifying space of the Lie group $G_2$ is 
$$H^*(BG_2;\mbbz)=\mbbz[y_4,y_{12}]\oplus\mbbz_2[y_6,y_{10}]$$ 
where the degrees are $|y_k|=k$.

\end{thm}

\section{Existence of Harvey-Lawson pairs}\label{hlpair}

Here we  illustrate an application of the topological 
results which we have proved in the previous section. 
By applying the Leray-Hirsch theorem  (e.g. \cite{sp}, \cite{hatcher}) 
Theorem~\ref {cupproductstructuregrassmannian} can be generalized from 
$G^{+}_{3}\R^7$ to the Grassmann bundle 
$\pi : \wt{M} \to M$. This is because of the fact that the Euler and 
Pontryagin classes $x=e(\nu)$, $y=p_{1}(\xi)$ are restrictions of the 
Euler and Pontryagin classes  $e(\V)$, $p_{1}(\Xi)$ of the corresponding  
universal bundles over $\wt{M}$ (cohomological extension property).

\beg{thm}\label{grassmannianbundle} $H^*(\wt{M};\Q)$ is an $H^{*}(M; 
\Q)$ module generated by $e(\V)$ and $p_{1}(\Xi)$. In other words the map 
$$a\otimes x + b\otimes y\mapsto  \pi^{*}(a) \cup e(\V) + \pi^{*}(b) \cup 
p_{1}(\Xi)$$ gives an isomorphism:
$$H^{*}(M )\otimes_{\Z} H^{*}G^{+}_{3}(\R^{7}) \longrightarrow H^*(\wt{M}).$$
\end{thm}

Next comes a corollary to the existence of Harvey-Lawson pairs. 
Recall that a manifold pair $Y^3\subset X^4\subset (M^7,\varphi)$ in a 
manifold with $G_2$-structure is called a {\em Harvey-Lawson pair} if 
the three form vanishes on the normal bundle of $X^4$ when restricted to 
$Y^3$. Note that these type of submanifolds are related to the Mirror-
duality  of \cite{as3} and \cite{as4}.

\beg{cor}\label{hl} Let $ X^{4}\hookrightarrow (M^{7},\varphi ) $ be any embedding of a closed smooth $4$-manifold into a manifold with $G_2$ structure  satisfying the property $\langle p_{1}(\nu X), [X] \rangle \neq \pm e[X]$, where $e[X]$ is the Euler characteristic. After a small isotopy of $X\subset M$, we can find  a nonempty closed smooth $3$-dimensional submanifold $Y^{3}\subset X^{4}$ such that $(X,Y)$ is a HL pair.
\end{cor}

\beg{proof} Consider the map 
$\Psi: Im(X, M)\times X \longrightarrow \wt{M}$ given by $\Psi(f,x)=f_{N}(x)$, where $Im(X,M)$ denotes the space of embeddings of $X$ into $M$. $f_N$ assigns the normal plane to $f(X)\subset M$ at the image of a point $x\in X$. By transversality (e.g. \cite{gp}) we can find a nearby isotopic copy $f$ of any embedding, such that $f_{N}:X\to \wt{M}$ is transverse to the submanifolds $\A SS_{-}\sqcup \A SS_{+} \sqcup \A SS_{0} $. Since 
$p_{1}(\nu X)\pm e(TX) \neq 0$, $f_{N}$ meets both $\A SS_{\pm}$ since their Poincar\'e duals are $ p_{1}(\Xi)\pm e(\V) $. Hence $Y^3:=f^{-1}(\A SS_{0}) \neq 0$, and by definition $(X, Y)$ is a HL-pair. 
\end{proof}

\begin{figure}[h] \bct \label{kavun}  
\includegraphics[width=.55\textwidth]{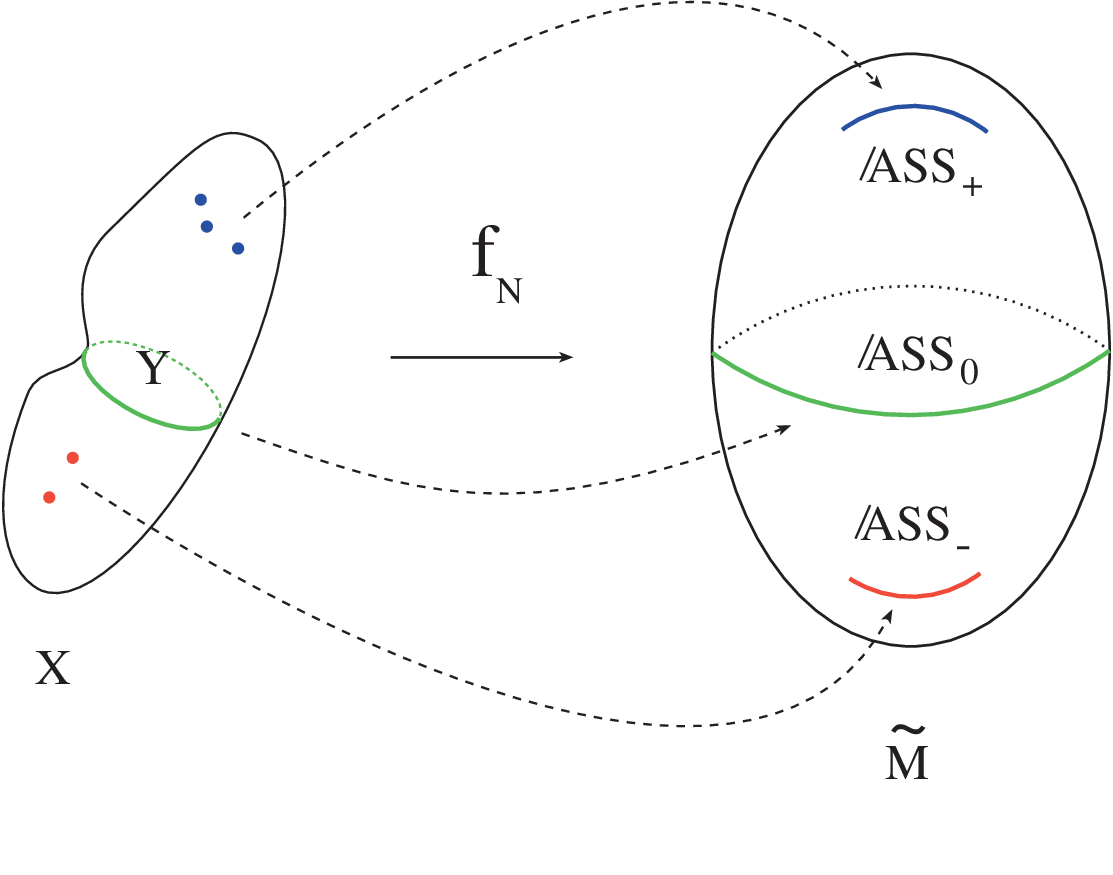}\\ \ect
  \caption{\, The map\, $\Phi: G^{+}_{3}\R^7 \to \R$ } 
\end{figure}

\beg{ex} By this corollary one can produce examples of HL pairs. Standard embeddings of $S^4$, $S^2\times S^2$ (after stabilization) into $\mbbr^7$ can be extended to an embedding into 
any manifold with $G_2$ structure via a coordinate diffeomorphism. Since these have 
trivial normal bundle and nontrivial Euler characteristic, they satisfy the hypothesis of the corollary and can be isotoped to a HL pair. 
Similarly an orientable closed surface $\Sigma_g$ of genus $g$ 
is embedded into $\mbbr^3$ with trivial normal bundle. So that  $\Sigma_g\times\Sigma_h$ can be embedded into any manifold with $G_2$ structure for $g\neq 1,h\neq 1$ and isotoped to a HL pair. 

\end{ex}

\vspace{.1in}

{\small 
\beg{flushleft} \textsc{Mathematics Department, Michigan State University, East Lansing, MI 48824
}\\
\textit{E-mail address:} \texttt{\textbf{akbulut@math.msu.edu},\textbf{kalafat@math.msu.edu}
}
\end{flushleft}
}




\end{document}